\documentstyle[12pt]{article}
\input{mssymb.sty}

\title{Quantization of Poisson pencils and generalized Lie algebras}
\author{D.Gurevich\\
Universit\'e de Valenciennes, \\ 59304 Valenciennes, France,\\
V.Rubtsov\\
Centre de Math\'ematiques, Ecole Polytechnique,\\
 91128 Palaiseau, France\\
and ITEP, Bol.Tcheremushkinskaya 25,\\
 117259 Moscow, Russia}
\date{}

\begin{document}

\newtheorem{proposition}{Proposition}
\newtheorem{definition}{Definition}
\newtheorem{remark}{Remark}
\newcommand{\gggg}{\frak g}
\newcommand{\bC}{\mbox{$\bf C$}}
\newcommand{\bR}{\mbox{$\bf R$}}
\newcommand {\co}{{\cal O}_x}
\newcommand {\cR}{{\cal R}}
\newcommand{\ug}{U(\gggg)}
\newcommand{\gug}{Gr\,U(\gggg)}
\newcommand{\us}{U(sl(2))}
\newcommand{\vv}{V^{\otimes 2}}
\newcommand{\uqs}{U_q(sl(2))}
\newcommand{\oI}{\overline{I}}
\newcommand{\uq}{U_q(\gggg)}
\newcommand{\sss}{SL(n)/S(L(n-1)\times L(1))}
\newcommand{\St}{\tilde S}
\newcommand{\sn}{sl(n)}
\newcommand{\ppb}{P_{\beta}}
\newcommand{\pqb}{P^q_{\beta}}
\newcommand{\wkm}{\wedge_-^k(V)}
\newcommand{\wm}{\wedge_-(V)}
\newcommand{\vb}{V_{\beta}}
\newcommand{\va}{V_{\alpha}}
\newcommand{\ww}{W^{\ot 2}}
\newcommand{\vaa}{V_{2 \alpha}}
\newcommand{\vva}{V_{\alpha}^{\ot 2}}
\newcommand{\wpp}{\wedge_+(V)}
\newcommand{\wkp}{\wedge_+^k(V)}
\newcommand{\wkpm}{\wedge_{\pm}^k(V)}
\newcommand{\wpm}{\wedge_{\pm}(V)}
\newcommand{\ot}{\otime}
\newcommand {\mod}{\mbox{$U_q(\gggg)$-Mod}}
\newcommand {\mos}{\mbox{$U_q(sl(2))$-Mod}}
\newcommand{\oi}{\mbox{$\overline{I}$}}
\newcommand{\aq}{A^q}
\def\ot{\otimes}
\def\de{\Delta}\def\SS{\frak S}
\newcommand{\ac}{A^q_c}
\def\sl{sl(2)}
 \def\qq{(q+q^{-1})}

\maketitle
\begin{abstract}

We describe two types of Poisson pencils generated by a linear
bracket and a quadratic one arising from a classical
R-matrix. A quantization scheme is discussed for each. The
quantum algebras are represented as the enveloping algebras of
``generalized Lie algebras''.
\end{abstract}

\section{Introduction}

It is well-known that Poisson pencils play an important role in the
 theory of integrable dynamical
systems. Recall that a {\em Poisson pencil} is a linear space of
Poisson brackets generated by two of them
\begin{equation}
\{\,\,,\,\,\}_{a,b}=a\{\,\,,\,\,\}_1+b\{\,\,,\,\,\}_2. \end{equation}
In this case the brackets $\{\,\,,\,\,\}_1$ and $\{\,\,,\,\,\}_2$ are called
{\em compatible}.

 However it was not clear what was a proper quantum analog of
Poisson pencil (1) in general. In the present paper we propose certain
 candidates for this role in a particular case
 when one of the brackets (say $\{\,\,,\,\,\}_1$)
 is linear
and the other ($\{\,\,,\,\,\}_2$) is quadratic
and defined by a classical ``canonical'' R-matrix (2).

Actually we consider two types of Poisson pencils (1) connected
with R-matrices
and discuss a way to quantize them.

Let us describe them but first we want to make a few remarks about our
notations.

All manifolds (or varieties)
will be considered over the field
$k=\bR$. By $Fun(M)$ we always mean the space of polynomials in the case
when $M=V$ is a linear space or the space of their restrictions to $M$
if $M$ is embedded as an algebraic
variety in $V$. In the last case, we say that a Poisson bracket is linear
(quadratic) if the Poisson bracket of the restrictions of two linear
 functions is the restriction of a linear (quadratic) function.

Note that all classical and quantum objects under consideration
have  natural holomorphic analogs over the field $k=\bC$.

{\bf Type 1.} Let $\gggg$ be a Lie algebra
and $\rho: \gggg\to Vect(M)$ be a representation of $\gggg$
in the space of
vector fields on a manifold $M$. Let us assign to
an element  $R\in \wedge^2(\gggg)$ a bi-vector field $\rho^{\ot 2}(R)$
and introduce a bracket
$$\{f,g\}_R=\mu<\rho^{\ot 2}(R),df\ot dg>.$$
Hereafter $\mu$ denotes the usual commutative multiplication in $Fun(M)$.
We are interested in the cases  when
this bracket is Poisson.

If it is,
 the manifold $M$ is equipped with another Poisson bracket $\{\,\,,\,\,\}$
and the elements of
$Im\,(\rho)$ preserve the latter bracket, i.e.,
$$\rho(X)\{\,\,,\,\,\}=\{\,\,,\,\,\}
(\rho(X)\ot id+id\ot \rho(X)),\,\forall X\in \gggg$$
 then we have a Poisson pencil (1) with
$\{\,\,,\,\,\}_1=\{\,\,,\,\,\}$ and $\{\,\,,\,\,\}_2=\{\,\,,\,\,\}_R$.

The  case under detailed examination is the following:
$\gggg$ is a
simple Lie algebra, $M=\co$ is the orbit in $\gggg^{*}$ of a highest weight
element $x\in\gggg^{*}$ (h.w. orbit for brevity) equipped
with the Kirillov-Kostant-Souriau bracket
$\{\,\,,\,\,\}=\{\,\,,\,\,\}_{KKS}$,
$R$ is assumed to be R-matrix (2) and $\rho=ad^{*}$.

{\bf Type 2.} Let $\gggg$ be a classical simple
Lie algebra (i.e., one of the following algebras
$sl(n),\, so(n),\, sp(n)$). Consider the Sklyanin-Drinfeld bracket
$\{\,\,,\,\,\}_{SD}$ (refered also as Poisson-Lie)
defined on the corresponding group $G$ as follows
$$\{f,g\}_{SD}=\{f,g\}_l-\{f,g\}_r,\,\,\{f,g\}_{r,l}=\mu<\rho_{r,l}^{\ot 2}(R),
df\otimes dg>,\,\, f,g\in Fun(G),$$
where $\rho_r (\rho_l)$ is the representation of the
 Lie algebra $\gggg$ in
the space $Fun(G)$ by the left- (right-) invariant vector fields
and $R$ is an R-matrix (2).

Let us consider an embedding of
 the group $G$ in the matrix
algebra $Mat(n)$. Then the bracket $\{\,\,,\,\,\}_{SD}$ can be extented
to the bracket
defined on the space $Fun(Mat(n))$ in a natural way but it is Poisson
iff $\gggg=sl(n)$. In this case the extended bracket
will be denoted $\{\,\,,\,\,\}_{2}$.

Note that the bracket $\{\,\,,\,\,\}_{2}$
is quadratic with respect to the basis given by the matrix elements.
Linearization of this bracket gives rise to a Poisson bracket
(denoted $\{\,\,,\,\,\}_1$) compatible with
the bracket $\{\,\,,\,\,\}_2$.
Thus we have a Poisson pencil defined in the space $Fun(Mat(n))$
and generated by the brackets $\{\,\,,\,\,\}_{1,2}$.

Finally we have two types of Poisson pencils arising from R-matrices.
The main difference between them is following. For the first type
Poisson pencil we start with a linear bracket and construct the quadratic
one by means of a classical R-matrix. In contrary, the linear bracket of
the second type Poisson pencil is appeared as a
linearization of a given quadratic bracket.
Unfortunately we do not know any scheme which could include both types.

Note that integrability of the
Toda chain is connected with the last type Poisson pencils.

Let us say a few words about the
quantization of Poisson pencils of both types.
Note that a quantization scheme of the first type Poisson pencils in the
case when $R$ is a classical (non-modified, cf. below) R-matrix
was proposed in \cite{GRZ}. The case when $R$ is
R-matrix (2) and $M$ are orbits in $\gggg^*$ of a special type was
examined in \cite{DG1}, \cite{DS2}.

The algebras quantizing the Poisson pencils under consideration
(quantum algebras for short) can be described as
follows. We describe the
quantum algebras corresponding to the brackets $\{\,\,,\,\,\}_2$
in terms of graded
 quadratic algebras\footnote{Let us precise
that we call a  {\em filtered quadratic algebra} (f.q.a.)
an algebra of the form
$T(V)/\{I\}$ where $V$ is a linear
space and $I\subset V^{\ot 2}\oplus V\oplus k$. Hereafter
we denote $T(V)$ the free tensor algebra of a linear space $V$
and $\{I\}$ denotes the ideal generated by a set $I$.
If $I\subset V^{\ot 2}$ we call the corresponding algebra
{\em graded quadratic} (g.q.a.).}.

For the second type Poisson pencils this
algebra is q-deformed symmetric algebra of the space $Mat(n)^{*}$.
For the first type ones
 it is a q-deformed g.q.a. which corresponds due to
Kostant's theorem (cf. \cite{LT}) to the h.w. orbit in $\gggg^{*}$.

To quantize the whole Poisson pencil (2)
we describe the correspoding quantum algebras in terms of f.q.a. which
look like the
enveloping algebras of the Lie algebras (with a parameter $h$ introduced
in the Lie bracket). Thus we have two parameter quantum algebras
$A_{h,q}$ (a parameter $q=e^{\nu}$ is one of ``braiding'' and
another parameter $h$ appears while  deforming a commutative
algebra of functions to an enveloping algebra).
The algebras $A_{h,q}$ will be called of the first or second type
accordingly to the type of the initial Poisson pencil.
(Let us emphasize that unlike our quantum algebras the famous quantum
groups $\uq$ are not  f.q.a.)

In this connection a natural question arises : whether it is
possible to represent the latter algebra itself as the
enveloping algebra of a q-deformed or ``generalized'' (in some sense)
Lie algebra structure
(in particular, what is a proper analog of Jacobi identity
in the general case)~?

There were number of attempts to find a proper q-analog of the
 Lie algebras.
In particular, objects of such a type have been considered
 in \cite{DOS} (cf. also the references therein). In \cite{W}
some generalized Lie structures were introduced in terms of
``quantum differential calculus''. A version of ``quantum'' and ``braided''
Lie algebras were proposed by Sh.Majid \cite{M}. However, in all these papers
the following question was ignored: whether the corresponding
enveloping algebra is a {\em flat}
 deformation\footnote{One says that a deformation
from $A_{0}$ to $A_{h}$ is flat if
$A_{h}/hA_{h}$ coincides with $A_{0}$ and
$A_{h}$ is a free $k[[h]]$-module. As for a deformation of g.q.a.
$A_0\to A_h$ it is flat iff $dim_{k[[h]]}\,A_h^p=dim\,A_0^p$
where $A_h^p$ is an homogeneous component of $A_h$ of
degree $p$. A deformation of a g.q.a. $A_{0}$ to f.q.a. $A_h$ is flat if
there exists a natural isomorphism of the PBW type between $Gr\,A_h$ and
$A_0$.}  of the initial
classical counterpart.

We propose our version of ``deformed Lie algebras'', called in the
paper ``generalized Lie algebras'', using as a criterium the flatness of
deformation of the
enveloping algebras. (In the present paper we consider only algebras
lying in {\em quasiclassical categories}, i.e. whose braiding operator is a
deformation of the flip\footnote{Note that
there exist non-quasiclassical categories, e.g, super-categories. Other
examples of non-quasiclassical rigid tensor or
quasitensor categories can be constructed by
means of non-quasiclassical solutions of the QYBE, introduced in
\cite{G2}, \cite{G3}.} though the below definition is valuable in the general
 case).
In this case we have two parameter family of associative algebras (the sense
of
the parameters is explained above) and as a quasiclassical limit of this family
we have a Poisson pencil.
Thus the existence of a Poisson pencil is a
necessary condition for the ``reasonable'' from the above point of view
deformation of a ordinary Lie algebra structure.
(Let us emphasize that this is not a deformation of one Lie algebra
to another and the usual deformation theory
can not be applied here.)

Thus, observing that the SD bracket (extended to $Fun(Mat(2))$) is not
compatible  with the linear bracket correspondind to the Lie algebra
$gl(2)$ it is easy to show that there does not exist any ``reasonable''
deformed
version of $gl(2)$ in the frames of the second type construction.
(As for the first type deformation of the Lie algebra $sl(2)$, cf. \cite{DG1}).

In any case if we want to get such a deformation of a Lie algebra
$\gggg$ that the deformation of the corresponding enveloping algebras
is flat we should restrict ourselves
to the h.w. orbit in $\gggg^{*}$ instead of the initial symmetric
algebra of the space $V=\gggg$ (for the type 1)
or take as the initial Lie algebra one given by the formula (5)
(for the type 2). In the latter case we can obtain a deformation of the
Lie algebra which forms a double Lie algebra with $gl(n)$.

The paper is organized as follows. In next Section we introduce two type of
 Poisson pencils
connected to R-matrix (2). In Section 3 we describe a scheme of their
quantization. The role
of ``generalized Lie algebras'' in the quantization procedure
is discussed in Section 4. An example of such a generalized Lie algebra
arising from a second type Poisson pencil is given in the last Section.

\section{Two types of Poisson pencils}

Recall that an element $R\in \wedge^2(\gggg)$ where $\gggg$ is a Lie algebra
is called a {\em classical R-matrix}
if the element
$$[[R,R]]=[R^{12}, R^{13}]+[R^{12}, R^{23}]+[R^{13}, R^{23}],$$
which always lies in $\wedge^3(\gggg)$, is equal to 0. We call $R$
 a {\em modified}  R-matrix if it satisfies the classical modified YBE,
i.e., the element $[[R,R]]$ is $\gggg$-invariant.

{\bf Type 1}.
Given a Poisson bracket $\{\,\,,\,\,\}$ on a manifold $M$, consider the space
$Vect(M,\,\{\,\,,\,\,\})$
of the vector fields $X\in Vect(M)$ preserving this bracket.
It is evident that they form a Lie algebra with respect to the
usual Lie bracket
$$X\ot Y \in Vect(M\,,\{\,\,,\,\,\})^{\ot 2}\to[X,Y]=X Y-Y X.$$

Let
$\rho:\gggg\to Vect(M,\,\{\,\,,\,\,\})$
 be a morphism of a finite dimensional Lie algebra $\gggg$
in the space $Vect(M\,,\{\,\,,\,\,\})$ equipped with the above Lie structure.
 Let us fix an element $R\in\wedge^2(\gggg)$.
Consider the bracket $\{\,\,,\,\,\}_R$ corresponding to this element.
If it is Poisson we have a Poisson pencil (1) since the brackets
$\{\,\,,\,\,\}$ and $\{\,\,,\,\,\}_R$ are always compatible.

It is, e.g., if $R$ is a classical non-modified R-matrix.
Let us assume now that $\gggg$ is a simple Lie algebra and $R$ is a
``canonical'' modified R-matrix
\begin{equation}
R=\sum_{\alpha \in \Delta_{+}}X_{\alpha}\wedge
X_{-\alpha} \in \wedge ^2 \gggg.
\end{equation}
Here $\{H_{\alpha},X_{\alpha},X_{-\alpha}\} $
is the Cartan--Weyl system  in the Chevalley normalization and $\Delta_{+}$
is the set of the positive roots of  $\gggg$ (other
normalizations of this R-matrix are possible as well).

Then the bracket $\{\,\,,\,\,\}_R$ is Poisson only
on certain manifolds (we say that such a manifold $M$ is of {\em
R-matrix type}).
All symmetric homogeneous $G$-spaces $M=G/H$ and the h.w. orbits
in any linear space possessing a $G$-module structure are of R-matrix type
(cf. \cite{DGM}, \cite{GP}).

Restricting ourselves to orbits $M=\co$
in $\gggg^*$ of R-matrix type we get a family of examples of Poisson pencils
taking as the bracket $\{\,\,,\,\,\}_1$
the Kirillov-Kostant-Souriau bracket $\{\,\,,\,\,\}_{KKS}$.

Other examples can be constructed as follows. Fix an even-dimensional
 linear space
$V$ equipped with a skew-symmetric non-degenerated pairing
$<\,\,,\,\,>:V^{\ot 2}\to k$ and
the Lie algebra $\gggg=sp(V)$ of linear transformations
preserving this form $(<Xu,v>+<u,Xv>=0)$. Introduce in the space
 Fun($V^*)$ a Poisson bracket $\{\,\,,\,\,\}_1$ in the following way~:
for linear functions $f,\,g\in V$ we set $\{f,g\}_{1}=<f,g>$.
Then taking as $R$ the R-matrix (2) corresponding to the Lie algebra $sp(V)$,
(which is isomorphic to the algebra  $sp(n)$),
we get a quadratic Poisson bracket  $\{\,\,,\,\,\}_R$ (it is Poisson since
$V\setminus 0$ is the h.w. orbit) and this bracket is compatible with the
 bracket $\{\,\,,\,\,\}_1$.

The following proposition is obvious.
\begin{proposition}
If $\rho(\gggg)$ falls into the space of (locally)
Hamiltonian vector fields then any
symplectic leave of the bracket $\{\,\,,\,\,\}_R$ lies in some
leave of the bracket $\{\,\,,\,\,\}_1$.
\end{proposition}

{\bf Type 2}.
Let us assume now that $G$ is a classical simple
 Lie group, i.e., one of the following
 Lie groups $SL(n),\, SO(n),\,Sp(n)$  and $R$ is R-matrix (2).
Let us consider on $G$ the SD bracket.

Note that this bracket can be represented in the following form
\begin{equation}
\{L_1\,\otimes,\, L_2\}_{SD}=
[R_{\rho},\,L_1\circ  L_2],\,\,L_1=L\otimes id,\,L_2=
id\otimes L.
\end{equation}
Here $L=(a_i^j),\,\, 1\leq i,j\leq
n=dim\, V$ is a matrix with the matrix elements $a_i^j$ (the low index
 labels rows),
$\rho : \gggg\to End\,(V)$ is the fundamental vector
 representation and $R_{\rho}=\rho^{\ot 2}(R)$.
 This bracket was originated by E.
Sklyanin namely in the form (3).

One can extend this bracket to the space $Fun(Mat(n))$ using the formula (3).
This means that the bracket $\{f,g\}_{SD}$ is
 well-defined for any
functions $f,g\in Fun(Mat(n))$. However the extended bracket is Poisson
iff $\gggg=sl(n)$. We denote the extended bracket $\{\,\,,\,\,\}_2$.

In what follows we consider in the frames of the second type construction
only the case $\gggg=sl(n)$.

Let us linearize now the bracket $\{\,\,,\,\,\}_2$ (the resulting bracket will
be
denoted $\{\,\,,\,\,\}_1$). The linearization can be defined as follows.
In the r.h.s. of the formula (3) we replace the matrix $L$
by $L+\lambda id$ (i.e., we substitute $a_{i}^{j}+\lambda
\delta_{i}^{j}$ instead
of $a_{i}^{j}$) and take only linear term in $\lambda$. In a compact form
the result
can be represented as follows
\begin{equation}
\{L_1\,,\, L_2\}_{1}=[R_{\rho},\,L_1+L_2].
\end{equation}
Then the latter bracket is
 Poisson. Moreover, the brackets
 $\{\,\,,\,\,\}_1$ and $\{\,\,,\,\,\}_2$ are  compatible.

This fact can be deduced from \cite{S} but it is easy to prove it directly
observing that the r.h.s. of (4) after the above substitution
does not contains any terms quadratic in $\lambda$.

Note that a linearizarion of the SD bracket on any group $G$ gives rise
according to the Drinfeld construction \cite{D}
to a Lie structure in $\gggg^{*}$
(which forms a so-called {\em Manin triple}
 toghether with the initial Lie algebra
$\gggg$).
So the corresponding (linear) Poisson-Lie bracket is well-defined in the
space $Fun(\gggg)$.

 As for the case when $\gggg=sl(n)$
the bracket (4) is an extention of the latter
bracket to the space $Fun(Mat(n))$.
In this case we have two Poisson
brackets $\{\,\,,\,\,\}_{1,2}$ defined
simultaneously in the space $Fun(Mat(n))$, (the first of them can be
reduced to $Fun(sl(n))$ and the second one to $Fun(SL(n))$.

\begin{remark} There exist (at least) two schemes to construct
 integrable dynamical systems. The first (Magri-Lenart scheme)
deals with bi-Hamiltonian
systems. (Recall that a dynamical system
is called bi-Hamiltonian if
it is Hamiltonian with respect to two different Poisson structures.)

Another scheme was developped by  Adler, Kostant and Symes.
In R-matrix interpretation, proposed by the
St-Petersbourg school, it deals with two Lie structures defined in the
same linear space (so-called double Lie algebras). We refer the
reader to \cite{S} for details.
\end{remark}

Note that for the second type Poisson pencils the symplectic leaves of
the first bracket  do not lie in the leaves of the second one. This  condition
is necessary for constructing the integrable dynamic systems.

Let us give now an explicit form of the second type Poisson brackets
$\{\,\,,\,\,\}_{1,2}$ under consideration.

The multiplication table for the second Poisson bracket is
$$
\{a_{k}^{i},a_{k}^{j}\}_{2}=a_{k}^{i}\,a_{k}^{j},\,\,
\{a_{i}^{k},a_{j}^{k}\}_{2}=a_{i}^{k}\,a_{j}^{k},\,\,i<j;
$$
$$
\{a_{i}^{l},a_{k}^{j}\}_{2}=0,\,\,
\{a_{i}^{j},a_{k}^{l}\}_{2}=2\,a_{i}^{l}\,a_{k}^{j},
\,\,i<k,\,j<l.
$$

It is not difficult to deduce from this the multiplication table for the first
Poisson bracket:
$$
\{a_{k}^{i},a_{k}^{j}\}_{1}=\delta_{k}^{i}\,a_{k}^{j}+a_{k}^{i}\,\delta_{k}^{j},\,\,
\{a_{i}^{k},a_{j}^{k}\}_{1}=\delta_{i}^{k}\,a_{j}^{k}+a_{i}^{k}\,\delta_{j}^{k},\,\,i<j;
$$
$$
\{a_{i}^{l},a_{k}^{j}\}_{1}=0,\,\,
\{a_{i}^{j},a_{k}^{l}\}_{1}=2\,(\delta_{i}^{l}\,a_{k}^{j}+a_{i}^{l}\,\delta_{k}^{j}),
\,\,i<k,\,j<l.
$$

It is useful to represent
the latter bracket under the following form
\begin{equation}
\{a_{i}^{j},a_{k}^{l}\}_{1}=\{{\bf R}(a_{i}^{j}), a_{k}^{l}\}_{gl}+
\{a_{i}^{j}, {\bf R}(a_{k}^{l})\}_{gl}.
\end{equation}
Here $\{\,\,,\,\,\}_{gl}$ is the Poisson-Lie
 bracket corresponding to the Lie algebra $gl(n)$
(namely, $\{a_{i}^{j},a_{k}^{l}\}_{gl}=
a_{i}^{l}\delta_{k}^{j}-a_{k}^{j}\delta_{i}^{l}$)
and ${\bf R}:W\to W$ is an operator defined in the space $W=Span(a_i^j)$ as
follows
${\bf R}(a_{i}^{j})=sign(j-i)\,a_{i}^{j}$ (we assume that $sign(0)=0$).

Actually we have equipped the space $W$
 with two Lie structures. The first one
is $gl(n)$ and the seconde corresponds to the Poisson bracket
$\{\,\,,\,\,\}_{1}$. Thus we have
a double Lie algebra in the space $W$ (our construction coincides with
one from \cite{S} up to a factor).

\section{Quantization scheme}

In the present section we discuss the quantization scheme for both
types of the Poisson pencils (2). This scheme consists of two steps. On
the first step we describe the quantum algebra for the bracket
$\{\,\,,\,\,\}_{2}$. This algebra can be represented as a g.q.a
for both types. However, there exists a difference between them.
If for the second type Poisson pencil
we deform the symmetric algebra of the initial space, for the first type
we deal with a deformation of the
g.q.a. coressponding to h.w. orbits in $\gggg^{*}$. Meanwhile for the
both cases the deformations are flat.

On the second step we are looking for linear terms converting the
mentioned g.q.a. into a f.q.a. in such a way that two properties are
fulfiled:\\
1. the deformation from the g.q.a. constructed before to f.q.a. is
flat;\\
2. the corresponding quasiclassical limit of this deformation coincides
with the given Poisson pencil.

{\bf Type 1}.
Let us fix a simple Lie algebra $\gggg$ and the corresponding quantum
group $\uq$.
Let us consider a linear space $V$ equipped with a $\gggg$-module structure
$\rho:\gggg\to End(V)$ and endow it with a structure of $\uq$-module
$\rho_{q}:\uq\to End(V)$. We assume that $\rho_{q}\to \rho$ when
$q$ goes to $1$.

Let us represent $\vv=\oplus \vb$ as a direct sum of the isotypical
$\gggg$-modules  ($\beta$ is the h.w. of the module $\vb$).
 If $\alpha$ is the h.w. of $V$ as $\gggg$-module then in
the above decomposition there is an irreducible $\gggg$-module $V_{2\alpha}$.
Let us set $\oI_{-}=\oplus_{\beta\not=2\alpha}\vb$ and
$\oI_{+}=V_{2\alpha}$.
Reproducing this procedure in the category of $\uq$-modules we can
introduce q-counterparts $\oI_{\pm}^{q}$ of the spaces $\oI_{\pm}$ (cf.
 \cite{DG2} for details).

Note that in virtue of the
Kostant theorem (cf. \cite{LT}) the function algebra on the h.w. orbit in
$V^{*}$ is just the algebra  $A_{0,1}=T(V)/\{\oI_{-}\}$.

As was shown in \cite{DS1} the deformation from the algebra
$A_{0,1}$ to its
q-deformed counterpart $A_{0,q}=T(V)/\{\oI_{-}^{q}\}$
is flat.

\begin{remark}
It is possible to define a q-analog of the  symmetric
algebra of the space $V$ in a similar way
but it is not, in general, a flat deformation of the initial object.
\end{remark}

It is not difficult to see that
the quasiclassical term of the deformation from $A_{0,1}$ to $A_{0,q}$
 coincides up to a factor depending on a normalization of the parameter
with the R-matrix bracket defined on the h.w.
orbit in $V^{*}$ (cf. also \cite{DG2}).
Thus the algebra $A_{0,q}$ is the quantum counterpart for the Poisson bracket
$\{\,\,,\,\,\}_R$ defined on the h.w. orbits in $V^{*}$.

In a particular case when  $V=\gggg$ ($\gggg$-module
structure of $\gggg$  is defined by ad-action) we are looking for a futher
(non-homogeneous) deformation of the algebra $A_{0,q}$. The latter deformation
will be realized
by means of the so-called {\em q-Lie bracket}.
In what follows we set $V=\gggg$.

Let $\alpha$ be the h.w. of $\gggg$ as $\gggg$-module.
Then in the decomposition $\vv=\oplus \vb$ where $V=\gggg$
there exists an unique irredicible
$\gggg$-module $\va$ belonging to the space $I_{-}$.

Let $\va^{q}$ be its analog in the category of $\uq$-modules
(i.e. the term in the decomposition of $\vv$ in the category of $\uq$-modules
which is a deformation of the $\gggg$-module $V_{\alpha}$, cf. \cite{DG2} for
details). Then we
introduce {\em q-Lie bracket} $[\,\,,\,\,]_{q}:\vv\to V$ by setting
$$
[\,\,,\,\,]_{q}|_{V^{q}_{\beta}}=0\,\, {\rm if}\,\,\vb^{q}\not=\va^{q}\,\,
 {\rm and}\,\,
[\,\,,\,\,]_{q}:\va^{q}\to V
$$
is a morphism in the category of $\uq$-modules (thus the bracket
is defined up to a factor).

Let us consider the algebra
$$
A_{h,q}=T(V)/\{\oplus_{\beta\not=2\alpha} Im(id-h[\,\,,\,\,]_{q})\vb^{q}\}.
$$

It is very plausible that the deformation from $A_{0,q}$ to $A_{h,q}$
is flat. Then the deformation
from $A_{0,1}$ to $A_{h,q}$ is flat as well.

Setting in the two parameter family
$h=a\hbar,\,\,q=e^{b\hbar}$ and looking for the Poisson bracket
corresponding
to this deformation one can reveal a bracket from the
pencil (1).
Thus, the family of algebras $A_{h,q}$
is the quantum object for the pencil (1) of the first
type defined on the h.w. orbit in $\gggg^{*}$.

{\bf Type 2}. Let us consider the base $\{a_{i}^{j}\}$ in the algebra
$Fun(Mat(n))$ consisting of the matrix elements.
We identify the spaces $Fun(Mat(n))$ and
$A_{0,1}=T(W)/\{I_{-}\}$ where $I_{\pm}$ is the subspace
of symmetric (skew) elements in $W^{\ot 2}$.
Recall that $W=Span(a^j_i)$.

A quantum (or q-) analog of this space can be described as follows.
Let us consider a solution $S=S_{V}:\vv\to\vv$ of the quantum Yang-Baxter
equation (QYBE)
$$
S^{12}S^{23}S^{12}=S^{23}S^{12}S^{23},\,\, S^{12}=S\ot id,\, S^{23}=id\ot S
$$
where $S=\sigma\rho^{\ot 2}({\cal R})$, ${\cal R}$ is a
quantum universal R-matrix
corresponding to the quantum group $\uq$, $\sigma$ is the flip and
$\rho:\uq\to End(V)$
is the vector fundamental representation.

Let us note that in the case under consideration ($\gggg=sl(n)$) the operator
 $S$ is of the ``Hecke type'', i.e., it possesses two eigenvalues.
 More precisely,
it has the following form
$$
S(a_i\otimes a_j)=(q-1)\delta_{i,j}a_i\otimes a_j+a_j\otimes a_i+
\sum_{i<j}(q-q^{-1})a_i\otimes a_j,
$$
where $\{a_i\}$ is a base in the space $V$.

There exists a natural way to introduce an operator
$S_{W}:W^{\otimes 2}\to W^{\otimes 2}$
which satisfies the QYBE
setting $S_{W}=S\ot {S^{*}}^{-1}$ where
$S^{*}:{V^{*}}^{\ot 2}\to {V^{*}}^{\ot 2}$
is defined with respect to the pairing
$$
<x\otimes y, a\otimes b> =<x,a><y,b>,\,\,x,y\in V,\,\,a,b\in V^{*}.
$$
Actually we treate the space $W$
as $V\otimes V^{*}$ and assume that the spaces $V$ and $V^{*}$ commut with
each other in the ordinary sense\footnote{Let us
remark that there exists another way
to extend the operator $S$ to the space $V\otimes V^{*}$ regarding
$V^{*}$ as
dual to $V$ in the rigid category generated by $V$ (cf. \cite{G3}).}.

In the base $\{a_i^j\}$ the operator $S_W$ is of the form
$$
S_{W}(a_{i}^{k}\ot a_{j}^{l})=S^{mn}_{ij}{S^{-1}}^{kl}_{pq}
(a_{m}^{p}\ot a_{n}^{q}) \,\,{\rm where}\,\,
S(a_{i}\ot a_{j})=S_{ij}^{kl}a_{k}\ot a_{l}.
$$

Let us point out an important property of the operator $S_W$: it
possesses 1 as an eigenvalue. This is the reason why it is natural to
introduce a deformed analog of the (skew)symmetric subspace of the space
$W^{\otimes 2}$, setting
$$I_{-}^{q}=Im(S_{W}-id),\,\,I_{+}^{q}=Ker(S_{W}-id),\,\,
A_{0,q}=T(W)/\{I_{-}^{q}\}.$$

The explicit form of these spaces is the following:
$$
I_-^q=Span(a_{k}^{i}a_{k}^{j}- qa_{k}^{j}a_{k}^{i},\,
a_{i}^{k}a_{j}^{k}-qa_{j}^{k}a_{i}^{k},\, i<j;
$$
$$
a_{i}^{l}a_{k}^{j}-a_{k}^{j}a_{i}^{l},\,
a_{i}^{j}a_{k}^{l}-a_{k}^{l}a_{i}^{j}-(q-q^{-1})a_{k}^{j}a_{i}^{l},\, i<k,j<l)
$$
and
$$
I_+^q=Span((a_i^k)^2,\,q a_{k}^{i}a_{k}^{j}+a_{k}^{j}a_{k}^{i},\,
qa_{i}^{k}a_{j}^{k}+a_{j}^{k}a_{i}^{k},\, i<j;
$$
$$
a_{i}^{j}a_{k}^{l}+a_{k}^{l}a_{i}^{j},\,
a_{i}^{l}a_{k}^{j}+a_{k}^{j}a_{i}^{l}+(q-q^{-1})a_{i}^{j}a_{k}^{l},\, i<k,j<l).
$$

\begin{proposition}(\cite{DS1}) The deformation from $A_{0,1}$ to
$A_{0,q}$ is flat.
\end{proposition}

Let us emphasize that considering the operators $S$ which correspond to
other classical simple algebras (and possess in this case 3 eigenvalues)
one can introduce in a similar way the spaces $I^q_{\pm}$. However
the deformation from $A_{0,1}$ to
$A_{0,q}=T(W)/\{I_{-}^{q}\}$ is not flat in this case.

It is easy to see that the Poisson bracket corresponding to the above
deformation coincides up to a factor with $\{\,\,,\,\,\}_{2}$. Thus
we can consider the algebra $A_{0,q}$ as a resulting object for the first
step of the quantization for the second type Poisson pencil.

As for a two-parameter family of associative algebras
$A_{h,q}$ quantizing the whole Poisson pencil (1) it is possible
 to introduce it by means of the mentioned above
substitution $a_i^j \to a_i^j+\lambda \delta _i^j$. Under this substitution
the space $I^q_-$ transforms to the space
$$
J^{h,q}=Span(a_{k}^{i}a_{k}^{j}- qa_{k}^{j}a_{k}^{i}-
h(\delta_{k}^{i}\,a_{k}^{j}+a_{k}^{i}\,\delta_{k}^{j}),$$
$$
a_{i}^{k}a_{j}^{k}-qa_{j}^{k}a_{i}^{k}-
h(\delta_{i}^{k}\,a_{j}^{k}+a_{i}^{k}\,\delta_{j}^{k}),\,\,i<j;
a_{i}^{l}a_{k}^{j}-a_{k}^{j}a_{i}^{l}=0,$$
$$
a_{i}^{j}a_{k}^{l}-a_{k}^{l}a_{i}^{j}-(q-q^{-1})a_{k}^{j}a_{i}^{l}-
hm\,(\delta_{i}^{l}\,a_{k}^{j}+a_{i}^{l}\,\delta_{k}^{j}),\, i<k,j<l),
$$
where $h=\lambda(q-1),\, m=1+q^{-1}$.

Then setting $A_{h,q}=T(W)/\{J^{h,q}\}$ we obtaine by construction
a flat deformation of the algebra $A_{0,q}$.

Thus we have quantized the both types of Poisson pencils
 under consideration and introduced the corresponding quantum algebras.
In next Section we discuss the
Lie algebra-like objects generating these quantum algebras.

\section{Generalized Lie algebras and flatness of deformation}

It is well-known that the axioms system of the ordinary Lie algebras
$\gggg$ allows one \\
1. to prove the PBW theorem
 and therefore to establish a flatness of deformation
from the symmetric algebra $Sym(\gggg)$ to the enveloping one $\ug$;\\
2. to introduce a notion of $\gggg$-module and
to construct a tensor category of $\gggg$-modules using the Hopf
structure of the algebra $\ug$.

In this Section we discuss a possible generalization of the notion of a
Lie algebra in the connection with flatness of deformations discussed above.

Let $V$ be a linear space and let $I$ be an arbitrary subspace in $\vv$.
Given a linear morphism $[\,\,,\,\,]:I\to V$, it is natural to introduce its
``enveloping algebra'' as follows
\begin{equation}
A=T(V)/\{Im(id-[\,\,,\,\,])I\}.
\end{equation}
(Let us note that in the classical case, i.e., when
$[\,\,,\,\,]$ is an ordinary Lie bracket,
our definition of the enveloping
algebra differs from the
usual one by the factor 2 in the bracket.)

 This  algebra is a filtered quadratic one. We consider
its adjoint graded algebra $Gr\,A$. It is evident that if there
 exists a graded isomorphism of the PBW type
$Gr\,A\to T(I)/\{I\}$ then the following relations are satisfied
\begin{equation}
([\,\, ,\,\,]^{12}-[\,\, ,\,\,]^{23})(I\otimes V\cap V\otimes I) \subset
I, \end{equation}
\begin{equation}
[\,\, ,\,\,]([\,\, ,\,\,]^{12}-[\,\, ,\,\,]^{23})(I\otimes V \cap V\otimes
I)=0. \end{equation}

If the g.q.a. $T(I)/\{I\}$ is Koszul (cf. \cite{PP})
then a reciprocal statement is true as well.

\begin{proposition} (\cite{PP}, \cite{BG})
Given a linear space $V$, a subspace $I\subset \vv$
and an operator (``Lie bracket'') $[\,\,,\,\,]:I\to V$. Then if
the g.q.a. $T(V)/\{I\}$ is Koszul and the
relations (7),  (8) are fulfilled then
there exists a graded isomorphism $Gr\,A\to T(V)/\{I\}$.
\end{proposition}

Let us point out that in this Proposition the bracket $[\,\,,\,\,]$ is
defined only on the subspace $I\subset \vv$. However in the cases under
consideration we always  have a pair of subspaces $I_{\pm}$ or $\oI_{\pm}$
(the space $I_{-}$ or $\oI_{-}$ plays the role of the space $I$ from
Proposition 4).
 This motivated us to introduce the following

\begin{definition}
 The data $(V,\,\vv=I_{+}\oplus I_{-},\,
[\,\,,\,\,]:\vv\to V)$ is called a generalized Lie algebra if the
following system of axioms is fulfilled\\
1.  $[\,\, ,\,\,]I_{+}=0$;\\
2. the relations (7), (8) take place with $I=I_{-}$;\\
3.  the g.q.a. $T(V)/\{I_{-}\}$ is Koszul.

Let ${\frak A}$ be a category. Then we say that  a generalized Lie
algebra is a
 ${\frak A}$-Lie algebra if one more axiom is fulfilled\\
4. the spaces $V,\,\, I_{\pm}$ are objects of the category ${\frak A}$
 and the bracket $[\,\, ,\,\,]$ is a morphism of this category.
\end{definition}

For the first time a definition of such a type was given in \cite{G4}.

The space $V$ equipped with a structure of a generalized
Lie algebra will be denoted by $\gggg$ and the enveloping algebra (6)
will be denoted by $\ug$.

Note that in the above definition we don't use any quantum R-matrix.
We formulate this definition in terms of two subspaces in the space $\vv$
but actually we use a quantum R-matrix to define the mentioned subspaces
in one or another way.

Note also that, since in frames of our construction (compared with one from
Proposition 3) the bracket $[\,\,,\,\,]$ is defined on the whole space $\vv$,
 we can introduce an analog of the adjoint operators by setting
$$
x\in V \to ad_{x}: V\to V\,\, {\rm where}\,\, ad_{x}y=[x,y].
$$
A way to extend such operator to a ``coadjoint vector field''
is discussed in \cite{DG2}.

\begin{proposition}
The algebras $A_{0,q}$ for generic $q$ are Koszul for the
both types of the above Poisson pencils.
\end{proposition}

It is so since the algebras $A_{0,1}$ are Koszul (for the fisrt type
algebras
this was proved in \cite{B} and for the second type ones it is well known)
and the fact that the deformation from $A_{0,1}$ to $A_{0,q}$ is flat
(cf. \cite{PP} for detail).

Thus the axiom 3 from the above definition is satisfied.

Using this fact we can state that the deformation from $A_{0,1}$ to
$A_{h,q}$ is flat iff the q-bracket satisfies the relation (7) and (8).
Assuming this, (for the case $\gggg=sl(2)$
this fact was proved in \cite{DG1}), we have by the construction
a ${\frak A}$-Lie algebra where ${\frak A}$ is the category of
$\uq$-modules.

As for generalized Lie algebras connected to the second type quantum algebras
$A_{h,q}$ it is natural to introduce
 the corresponding bracket by setting $[\,\,,\,\,]I_+=0$ and
$$
[\,\,,\,\,](a_{k}^{i}a_{k}^{j}- qa_{k}^{j}a_{k}^{i})=h(\delta_{k}^{i}a_{k}^{j}+
\delta_{k}^{j}a_{k}^{i}),\,
[\,\,,\,\,](a_{i}^{k}a_{j}^{k}-qa_{j}^{k}a_{i}^{k})=h(\delta_{j}^{k}a_{i}^{k}+
\delta_{i}^{k}a_{j}^{k}),\,\, i<j;
$$
$$
[\,\,,\,\,](a_{i}^{l}a_{k}^{j}-a_{k}^{j}a_{i}^{l})=0,\,
[\,\,,\,\,](a_{i}^{j}a_{k}^{l}-a_{k}^{l}a_{i}^{j}-(q-q^{-1})a_{i}^{l}a_{k}^{j})=$$
$$
m h(\delta_{k}^{j}a_{i}^{l}+
\delta_{i}^{l}a_{k}^{j}),\,i<k,j<l.
$$

It is evident by the construction that the relation (7) and (8) are satisfied.
Thus the latter bracket defines in the space $W$ a structure of a generalized
 Lie algebra indeed.

Note that the enveloping algebra of a generalized Lie algebra
does not have in general  any
``generalized bialgebra structure'' and therefore it is not
 possible to define a tensor
product of two modules over this associative algebra
in the category of $\ug$-modules.
Let us precise ungoing into details
that under ``generalized bialgebra structure'' we
mean a bialgebra structure which is defined in a tensor or quasitensor
category. Certain algebras of such a type
(called braided groups) were constructed by Sh. Majid but his braided
groups arising from the quantum group $\uq$ are not deformational objects
(in particular, it is not possible to define the corresponding Poisson
bracket).

However there exists a class  of  generalized Lie algebras
whose enveloping algebras admit ``generalized bialgebra structures'',
namely   S-Lie algebras (introduced in \cite{G1} under the name of generalized
Lie algebras).
 Their enveloping algebras
are objects of the symmetric tensor categories (for these categories
the ``commutativity morphism'' $S$ is involutive: $S^2=id$).
We refer the reader to \cite{GRZ} for a definition of  $S$-Lie algebras.

Completing this Section we want to compare our definition of a
generalized Lie algebra with close notions introduced by other authors.

Let us first note that construction of
the ``quantum Lie algebras'' from \cite{M}
is based on the following fact. There exists a subspace $L$ in the
quantum group $\uq$ closed with respect to the quantum
adjoint action of $\uq$ on itself.
Then a quantum Lie structure on $L$ is defined by the restriction of this
 action on $L$. A similar construction is given in \cite{DOS}
(cf. also the references therein) where a definition of the enveloping
algebra of such type algebras  is discussed as well. However the
latter algebra is not a flat deformation of any classical counterpart.

The same feature is inherent in the so-called ``braided Lie algebras'' from
\cite{M} (whose enveloping algebras are discussed in \cite{M} as well).
The definition of such a Lie algebra-like object consists in a linear
space beloning to a braided category and equipped with
three operators: a Lie bracket, a coassociative coproduct
and a counit which satisfy some axioms.

In contrary, we use only an operator of a Lie bracket in our definition
but we impose the axiom of Koszulity to ensure flatness of the deformation
of the enveloping algebra.

Let us emphasize that the most controversial question in all attempts
to introduce Lie algebras-like structures is a ``reasonable'' form
of the Jacobi identity.

It is worth to note that in the general case we can not reduce the difference
between two brackets
$[\,\,,\,\,]^{12}-[\,\,,\,\,]^{23}$  in the relations (7)-(8) to one of them.
It is possible only for $S$-Lie algebras whose  Jacobi identity has
more familiar form
$$[\,\,,\,\,][\,\,,\,\,]^{12}(id+S^{12}S^{23}+S^{23}S^{12})=0$$
or
$$[\,\,,\,\,][\,\,,\,\,]^{12}= [\,\,,\,\,][\,\,,\,\,]^{23}(id-S^{12}).$$

Under the latter form the Jacobi identity has appeared in the paper
\cite{W} in the frames of the quantum differential calculus.
However the quantum differential calculus of Woronowicz is not a result of any
flat
 deformation of the ordinary one and therefore the latter form
of Jacobi identity is not motivated by the deformation theory
(although the basic object of the quantum differential calculus,
 the algebras of quantized function $Fun_q(G)$, has a deformational
nature and  arises as a quantization of a Poisson bracket).

\section{Example: $n=2$ (type 2)}

Consider the algebra $Mat(2)$ consisting of the matrices
$\left(\begin{array}
{cc} a&b\\c&d\end{array}\right)$. The space $Fun(Mat(2))$ is generated by
the matrix elements $a,\,b,\,c,\,d$. The bracket $\{\,\,,\,\,\}_{2}$
in this case is determined by the following multiplication table
$$
\{a,b\}_{2}=ab,\,\,\{a,c\}_{2}=ac,\,\,\{a,d\}_{2}=2bc,
$$
$$
\{b,c\}_{2}=0,\,\,\{b,d\}_{2}=bd,\,\,\{c,d\}_{2}=cd.
$$

The linearization of this bracket gives rise to the linear bracket
$\{\,\,,\,\,\}_{1}$~:
$$
\{a,b\}_{1}=b,\,\,\{a,c\}_{1}=c,\,\,\{a,d\}_{1}=0,
$$
$$
\{b,c\}_{1}=0,\,\,\{b,d\}_{1}=b,\,\,\{c,d\}_{1}=c.
$$

Let us describe now the corresponding quantum algebras.
First we consider the space $W=Span(a,\,b,\,c,\,d)$ equipped with the
operator
$$S_{W}= S\otimes S^{-1}:W^{\otimes 2}\to
W^{\otimes 2}.$$
where the solution $S$ of the QYBE has the following form
$$
\left(\begin{array}
{cccc}
q&0&0&0\\0&q-q^{-1}&1&0\\0&1&0&0\\0&0&0&q\end{array}\right).
$$

Let us represent now the space $W^{\otimes 2 }$ as a direct sum of two
subspaces
$$
I_{-}^{q}=Span(ab-qba,\,ac-qca,\,
bd-qdb,\,\\cd-qdc,\,bc-cb,\,ad-da-(q-q^{-1})cb),
$$
$$
I_{+}^{q}=Span(a^{2},\,b^{2},\,c^{2},\,d^{2},\,qab+ba,\, qac+ca,\, qbd+db,
$$
$$
qcd+dc,\,ad+da,\, bc+cb+(q-q^{-1})ad),
$$
The algebra $A_{0,q}=T(W)/\{I_{-}^{q}\}$ is a deformation of the symmetric
algebra  of the space $W$ ``in direction'' of the bracket
$\{\,\,,\,\,\}_{2}$.
Now we deform the latter algebra in the second time
introducing the algebra $A_{h,q}$ as the quotient
algebra of $T(V)$ over the ideal generated by the elements
$$
ab-qba-hb,\,ac-qca-hc,\, bd-qdb-hb,
$$
$$
cd-qdc-hc,\,bc-cb,\,ad-da-(q-q^{-1})cb.
$$

Let us describe the space $I_{-}^{q}\otimes W \cap W\otimes
I_{-}^{q}$ explicitly.
It is generated by the following four
elements
$$
(ab-qba)c-(ac-qca)b+q^2(bc-cb)a=a(bc-cb)-b(ac-qca)+qc(ab-qba),
$$
$$
(ab-qba)d-q(ad-da-(q-q^{-1})cb)b+q(bd-qdb)a+(q^{2}-1)(bc-cb)b=
$$
$$
a(bd-qdb)-qb(ad-da-(q-q^{-1})cb)+qd(ab-qba),
$$
$$
(ac-qca)d-q(ad-da-(q-q^{-1})cb)c+q(cd-qdc)a=
$$
$$
a(cd-qdc)-qc(ad-da-(q-q^{-1})cb)
+qd(ac-qca)+(q^{2}-1)c(bc-cb),
$$
$$
(bc-cb)d-q(bd-qdb)c+q(cd-qdc)b=b(cd-qdc)-c(bd-qdb)+q^2d(bc-cb).
$$

It is instructive to check by direct computations that the relations (7) and
(8)
are satisfied.
Let us introduce now a generalized Lie bracket accordingly to the above scheme
by setting
$$
[\,\,,\,\,]I_{+}=0,\,[\,\,,\,\,](qab-ba)=hb,...,
[\,\,,\,\,](ad-da+(q-q^{-1})cb)=0.
$$
{}From this one can immediately deduce the following multiplication table
for the bracket $[\,\,,\,\,]$:
$$
[a,a]=[b,b]=[c,c]=[d,d]=[a,d]=[d,a]=[b,c]=[c,b]=0,\,[a,b]=[b,d]=Mb,
$$
$$
[b,a]=[b,d]=-Mqb,\,[a,c]=[c,d]=Mc,\,[ca]=[d,c]=-Mqc,
$$
with $M=h(1+q^{2})$.

The reader can readily compare this table with one for the ``braided Lie
algebra sl(2)'' from \cite{DG1}.

Completing the paper we would like to stress once more
that the second type Poisson bracket
$\{\,\,,\,\,\}_2$
 is not compatible with one $\{\,\,,\,\,\}_{gl}$. The reader can easily
to check this for the case $n=2$ observing that
$$\{a,\{b,d\}_{gl}\}_2+\{a,\{b,d\}_2\}_{gl}+c.p.\not=0$$
(general case follows from this one).
Therefore any reasonable (in the above sense)
deformation of Lie algebra $gl(n)$ does not exist in the frame of the second
type construction.

{\bf Acknowledgements.} The authors want to thank Mathematical
Departement of Universit\'e Claude Bernard, Lyon 1, where they begin to work
on the paper. They are also acknowledged the hospitality of Centre de
 Math\'ematique de l'Ecole
Polytechnique where the final version of the paper was prepared.
V.R. was partially supported by RFFR-MF

\end{document}